\documentclass[12pt]{article}

\usepackage[a4paper]{geometry}
\usepackage{ajc-FR}
\usepackage{amsmath,amsfonts,latexsym}
\usepackage{graphicx}
\usepackage{hyperref}
\usepackage{placeins}

\Volume{XY(Z)}
\Year{2022}
\firstpage{1}
\lastpage{9}
\Revised{29 November 2026}
\runninghead{\copyright \, F. Rowley November 2026\,\, -- ILBMRNUSS -- 1.05
DRAFT arXiv version}

\numberwithin{equation}{section}

\begin{document}

\title{\bf Improved Lower Bounds for\\Multicolour Ramsey Numbers\\using SAT-Solvers}
\author{Fred Rowley\thanks{formerly of Lincoln College, Oxford, UK.}}
\Addr{West Pennant Hills, \\NSW, Australia.\\{\tt fred.rowley@ozemail.com.au}}

\date{\dateline{17 September 2022}{DD Mmmmm CCYY}}

\maketitle


\begin{abstract}
This paper sets out the results of a range of searches for linear and cyclic graph colourings with specific Ramsey properties.  The new graphs comprise mainly 'template graphs' which can be used in a construction described by the current author in 2021 to build linear or cyclic compound graphs with inherited Ramsey properties.  These graphs result in improved lower bounds for a wide range of multicolour Ramsey numbers.

Searches were carried out using relatively simple programs (written in the language `C') to generate clauses for input to the PeneLoPe and Plingeling parallel SAT-solvers.  When solutions were found, the output from the solvers specified the desired graph colourings. 

The majority of the graphs produced by this work are `template graphs' with parameters in the form $(k,k,3)$ or $(k,l,3)$ with $k \ne l$.  Using these template graphs in familiar constructions, it has been possible to demonstrate significant improvements for lower bounds for most $R_r(k)$ for $5 \le k \le 9$ and $r \ge 4$.  These improvements provide correspondingly increased lower bounds on $\Gamma(k) = \lim_{\substack{r \rightarrow \infty}} R{_r}(k)^{1/r}$.  

We also show that $R_3(8) \ge 6534$ \footnote{Correction November 2026.  The author apologises to readers who may have seen previous versions of this paper with an incorrect value here.} and $R_3(9) \ge 15041$.  Other new lower bounds include $R(3,6,6) \ge 338$ and $R(3,8,8) \ge 941$, based on non-template cyclic graphs, and the interesting particular cases $R(3,4,5,5) \ge 729$ and $R(3,5,5,5) \ge 1429$. 

A spreadsheet containing specimens of many of the graphs mentioned here will be attached as an ArXiv ancillary file.  


\end{abstract}


\small \copyright Fred Rowley,  September 2022.
\bigskip

\section{Introduction}
This paper records the results of a series of automated constructions of multicolour classical Ramsey numbers, seeking to produce improved lower bounds on a range of diagonal and non-diagonal numbers. The construction of improved 'template graphs' leads to improvements in lower bounds on $\Gamma(k) = \lim_{\substack{r \rightarrow \infty}} R{_r}(k)^{1/r}$ for some small values of $k$

Although SAT-solvers were used to produce the results set out here, no knowledge of them or their use is necessary to understanding this paper.  Readers interested in the particular SAT-solvers used may wish to refer to the URL's mentioned in the Acknowledgements section, for more information.  

The paper is not intended to be fully standalone, and readers are advised to refer to \cite{Rowley3} where necessary detail is lacking here. However, we will repeat some terminology and definitions for convenience in Section~2.  

A brief history of the most relevant linear graph constructions is included in Section~3.  The results are summarised in Section~4, and set out in detail in Sections 5, 6 and~7.

A few brief conclusions are included in Section 8.

\section{Definitions and Notation}

Readers familiar with the terminology of this subject will no doubt omit a full reading of this section, and refer back if necessary.  Even so, we note that in this paper, 

$K_n$ denotes the complete graph of order $n$.  

If $U$ denotes a complete graph $K_m$ with $m$ vertices $\{u_0, {\dots} , u_{m-1}\}$, then a {\it (q-)colouring} of $U$ is a mapping of the edges $({u_i}, {u_j})$ of $U$ into the set of integers $s$ where $1\le s \le q$.  

The {\it length} of the edge $({u_i}, {u_j})$ is defined as $\mid j - i \mid$.  A length is often referred to as an {\it edge-length} in this paper, for clarity.

A colouring of $U$ is {\it linear} if and only if the colour of any edge $({u_i}, {u_j})$ depends only on the length of that edge.  In such a case the colour of an edge of length $l$ may be written $c(l)$, or $c(l,U)$ where necessary to avoid ambiguity.  

A colouring of $U$ is {\it cyclic} if and only if (a) it is linear, and (b) $c(l)=c(m-l)$ for all $l$ such that $1 \le l \le m-1$.  


The {\it clique number} of graph $U$ in colour $s$ is the largest integer $i$ such that $U$ contains a subgraph which is a copy of $K_i$ in that colour.  

We define a {\it Ramsey graph} $U(k_1, {\dots} \, ,k_r; m)$, with all $k_s \ge 2$, as a complete graph of order $m$ with an edge-colouring such that for each colour $s$, where $1 \leq s \leq r$, there exists no monochrome copy of a complete graph $K_t$ which is a subgraph of $U$ in the colour $s$ for any $t \geq k_s$.  Equivalently, the maximum order of any such copy of $K_t$ in any colour $s$ is strictly less than $k_s$.  Such a graph $U$ may conveniently be described as a $(k_1, {\dots} , k_r; m)${\it-graph}.  

The {\it Ramsey number} $R(k_1, \dots \, , k_r)$ is the unique lowest integer $m$ such that no \newline
 $U(k_1, {\dots} \, , k_r; m)$ exists, and its existence is proved by Ramsey's Theorem.  When all the $k_i$ are equal, this may be written $R_r(k)$ and is referred to as a {\it diagonal} Ramsey number.  

If the subset comprising all the edge-lengths of a linear Ramsey graph $U(k_1, {\dots} \, ,k_p, 3; m)$ with colour $p+1$ contains the edge-length $m-1$, we may call it a {\it (triangle-free or tf-) template}  for $U$. It should be clear that this definition ensures the subset is triangle-free. The subset and the corresponding induced subgraph of $U$ can be identified with each other without confusion.  
 
The construction described in \cite{Rowley3} (2021) allows the creation of linear Ramsey graphs with specific properties by combining the attributes of two linear graphs, which we will refer to as {\it prototypes}.  We may refer to graphs constructed in this way generically as {\it compound graphs}.  The approach in that paper involves the use of two prototype graphs, one of which has specific properties and is referred to as a {\it template graph}.  A template graph contains a triangle-free template in one colour (as defined above) but also has the property that it can be used in the 2021 construction to produce indefinitely large compound graphs, without increasing clique numbers in the other colours.  Much more detail about that construction can be seen in the paper, but is omitted here.

\section{Very Brief History}

In \cite{Giraud1}, Giraud effectively established that various infinite series of cyclic Ramsey graphs could be constructed, adding one colour at a time, and that the orders of those graphs could be determined by a simple formula based on the number of colours.  

Among other things, his results set lower bounds for a wide range of multicolour cyclic Ramsey numbers $R_r(k)$ of unlimited size, and for any $k$.  If we write $L^C(k_1, k_2, \dots , k_p)$ as the largest m such that a cyclic edge-colouring of $K_m$ in $p$ colours exists, which contains no instances of $K_{k_i}$ in any colour $i$, where $1 \le i \le p$, then Giraud proved that 

\begin{equation}
	L^C(k_1, k_2, \dots , k_p, k_{p+1}) \ge (2k_{p+1}-3).L^C(k_1, k_2, \dots , k_p)
	\label{f100}
\end{equation}

These bounds were improved later by Abbott and Hanson in \cite{AbbH}, by showing that two linear multicolour graphs could be compounded in a more effective way, where all of the values of $k_i$ are identical (across both graphs).  Accordingly, if we write $L^L_p(k)$ as the largest m such that a linear edge-colouring of $K_m$ in $p$ colours exists, which contains no instances of $K_k$ in any colour, then the formula in \cite{AbbH} showed that

\begin{equation}
	L^L_{(p+q)}(k) \ge ((2.L^L_p(k) - 1).(2.L^L_q(k) - 1)+1)/2
	\label{f101}
\end{equation}

A paper by the current author in 2017 \cite{Rowley1} removed the need for any condition of equality of the $k_i$ in the general linear case.  Thus, if we now write $L^L(k_1, k_2, \dots , k_p)$ as the largest $m$ such that a linear edge-colouring of $K_m$ in $p$ colours exists, which contains no instances of $K_{k_i}$ in any colour $i$, then the formula in \cite{Rowley1} demonstrated that:


\begin{equation}
	L^L(k_1, k_2, \dots , k_{p+q}) \ge ((2.L^L(k_1, k_2, \dots , k_p) - 1).(2.L^L(k_1, k_2, \dots , k_q) - 1)+1)/2
	\label{f102}
\end{equation}

and in the cyclic case that:


\begin{equation}
	L^C(k_1, k_2, \dots , k_{p+q}) \ge ((2.L^C(k_1, k_2, \dots , k_p) - 1).(2.L^C(k_1, k_2, \dots , k_q) - 1)+1)/2
	\label{f103}
\end{equation}

The two formulae above include the results of \cite{Giraud1} and \cite{AbbH} as special cases.

In 2021, a further construction was introduced (in \cite{Rowley3}), which increases the size of the compound graphs produced in a large proportion of cases.  This construction involves the use of a graph colouring featuring a specific form of tf-template (defined above).  

We can write $L^T(k_1, k_2, \dots , k_p, 3)$ as the largest m such that a linear edge-colouring of $K_m$ in $p+1$ colours exists, which produces no instances of $K_{k_i}$ in any colour $i$ under the conditions of multiple repetition specified in \cite{Rowley3},  and where $k_{p+1} = 3$ and $(p+1)$ is the colour of the tf-template.  We could equally say that $L^T(k_1, k_2, \dots , k_p, 3)$ is the order of the largest template graph with clique numbers strictly below the specified parameters.  

The formula in \cite{Rowley3} demonstrates that


\begin{equation}
	L^L(k_1, k_2, \dots , k_{p+q}) \ge (L^T(k_1, k_2, \dots , k_p, 3) - 1).(L^L(k_1, k_2, \dots , k_q) - 1) + 1 + \phi(T)
	\label{f104}
\end{equation}

The function $\phi(T)$ is dependent on the template concerned, being one less than the lowest edge-length coloured in the template colour.  

The construction in \cite{Rowley1} demonstrates that given any Ramsey graph $G(k_1, k_2, \dots , k_p)$ of order $m$ it is always possible to produce an effective template graph of order $(2m -1)$.  The 2021 formula therefore never produces a lower outcome than the 2017 formula.  

These approaches can be used repeatedly, producing in each case an infinite series of graphs with known orders.  In the diagonal case, the growth rate of the orders of the graphs in such a series can often be seen to approach a fixed limit as the number of colours tends to infinity, which we will refer to as the {\it limiting growth rate} of the series.  Where the series arises from the repeated use of the same construction with the same prototype graph, we may associate the limiting growth rate with that combination.  

\section{Overview of Current Results}

Since the author's 2021 construction requires the use of template graphs with specific properties, the current work has sought to increase the size of available template graphs, with some useful results with results set out in Sections 5, 6 and 7 below.  

It is noted that all graphs and templates mentioned in this paper have been tested, with regard to both their clique numbers and their usefulness as templates, independently of the process used to generate them, except in the `quadrupled' cases in Section~7.  

The range of searches carried out was partly constrained by the number of clauses necessary to define the problem to the SAT-solver, and the time taken to generate the files containing them.  For that reason, not all possibilities have been examined, and we are not yet able to say definitively whether any of the template graphs is of maximal order.  However, it appears that (e.g.) there is no template graph with parameters $(4,4,3)$ that outperforms the simpler approach represented by formula \ref{f102} or \ref{f103}.

In many cases, specimen template graphs included in the attached spreadsheet file exhibit the symmetries needed to produce cyclic compound graphs (provided that the other prototype is cyclic).

Using these template graphs, we can establish improved lower bounds on a wide range of `small' diagonal multicolour Ramsey numbers and the limits $\Gamma(k)$.  Specific new results are detailed in Section 5, including some improved lower bounds for $\Gamma(k)$.  These were derived directly from the parameters of the template graph used, as described in \cite{Rowley3}.  

In Section 6, one of the key Tables from \cite{Rowley3} is updated to reflect the improved results in particular cases for the previous construction.  

In Section 7, we deduce new lower bounds for $R_3(8)$ and $R_3(9)$.  

\section{Some New Templates and Improved Lower Bounds}

We will set out below the parameters of several small template graphs as described above.

For smaller graphs, (broadly where all parameters $k_i < 6$) the SAT-solver was simply allowed to specify all the constraints necessary to define such a graph, and colours at all distances were freely variable.  Two smaller examples are shown in Figure 1 below.  However, it was found that this approach did not work effectively for graphs of order approaching 100 using three or more colours.  

\begin{figure}[!ht] 
  \begin{center}
    \includegraphics{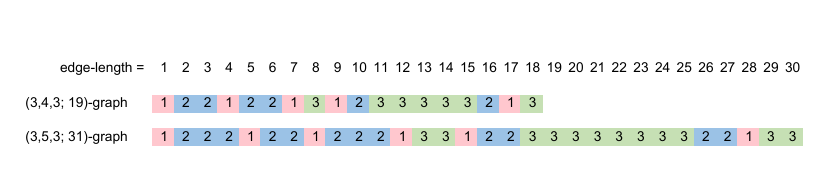}
  \end{center}
  \caption{\label{Fig:01} Two small template graphs.}
\end{figure}

For larger graphs, a different approach was adopted.  This involved first choosing a known cyclic graph of order $n$ as a prototype, and then colouring the edges in the interval $[1, n-1]$ in the same pattern as the prototype, according to their length.  Length $n$ was coloured in the template (special) colour.  

A variable $t$ was selected as the basis for each individual search based on the chosen prototype.  The order the template graph was determined as $2n+t$.   

Edge-lengths in the interval $[n+1,n+t]$ were allowed to vary (subject to the constraints) during the search.

Edge-lengths in the interval $[n+t+1, 2n-1]$ were coloured in the template colour.  The colouring for lengths greater than $n-1$ was constrained to be symmetrical, by reflection.  As a result, any lengths in the interval $[2n, 2n+t-1]$ were also variable, although it was not necessary to compute their colours separately in practice.  

To keep file sizes manageable, the constraints programmed under the approach above were initially specified simply to produce cyclic graphs with the desired parameters, and then a second program added constraints designed to allow the pattern of the colouring to repeat in the manner required for a template file.  This sometimes required a number of iterations, in order to specify these repetition constraints completely.

Specifications for many of the graph colourings constructed in this way are included in the ancillary file linked to this paper.  Table 1 below outlines some of the most significant examples, and some of the new lower bounds they demonstrate.  
 
\begin{table}[!ht]
  \begin{center}
    \includegraphics{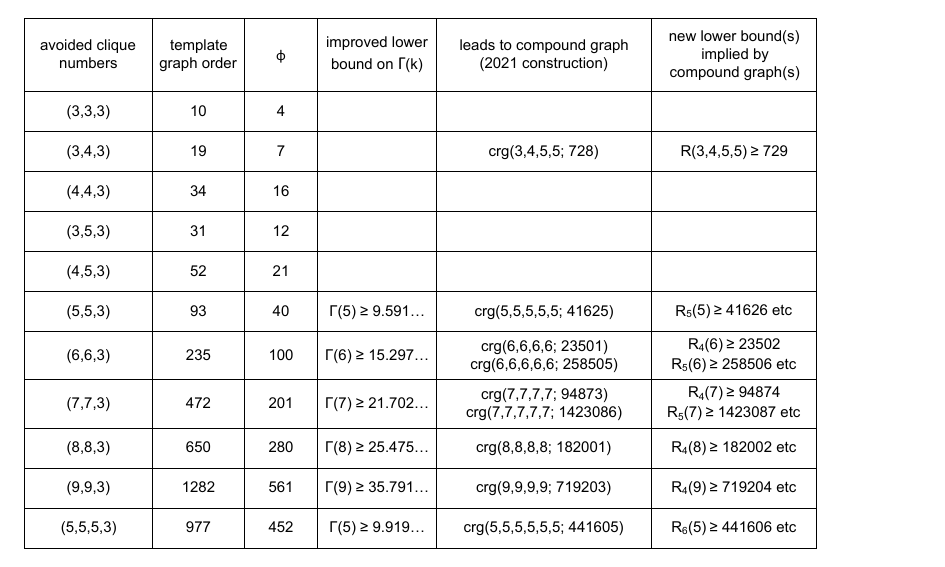}
  \end{center}
  \caption{\label{Table:01} Template graphs and derived lower bounds.}
\end{table}

\FloatBarrier

In Table 1, the abbreviation "crg" indicates a cyclic Ramsey graph with the specified parameters.  The template graph of order 977 was derived by extending a $(5,5,5; 453)$-graph derived by Exoo, and mentioned in \cite{RadzDS}.  The template graph of order 93 with parameters $(5,5,3)$ was generated by extending a cyclic $(5,5; 41)$-graph. The template graph of order 235 was derived by extending the $(6,6; 101)$-graph described by Kalbfleisch in \cite{Kalb}.  Other graphs for larger $k$ were similarly based on Paley graphs, using the Mathon-Shearer `doubling' process where appropriate.  

Table 2 below shows the new lower bounds obtained in some diagonal cases by using the new template graphs in the 2021 construction from \cite{Rowley3}.  That construction in its simplest form can also produce cyclic graphs with parameters $(5,5,5,3; 1429)$, $(6,6,3; 335)$, $(7,7,3; 673)$, $(8,8,3; 930)$ and $(9,9,3; 1843)$.  

We note here that the same searches also produced cyclic non-template graphs of orders 337 and 940, showing that $R(3,6,6) \ge 338$ and $R(3,8,8) \ge 941$.


\section{Larger Graph Orders and Lower Bounds on $\Gamma(k)$}

The tables below update the results of the previous paper \cite{Rowley3}, allowing for the inclusion of several new graphs constructed using the method mentioned above.  Only diagonal graphs are dealt with here.  

The table shows the best known linear graph orders in selected cases.  Bold green text indicates numbers exceeding those shown in the 2021 Radziszowski Dynamic Survey \cite{RadzDS}. 

\begin{table}[!ht]
  \begin{center}
    \includegraphics{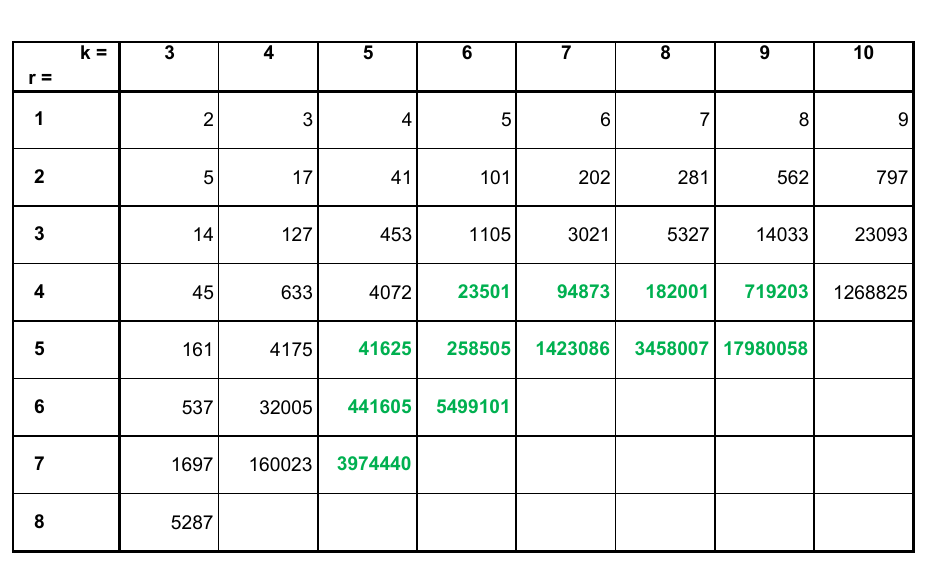}
  \end{center}
  \caption{\label{Table:02} Highest order of linear Ramsey graphs $U_r(k)$ known to the author.}
\end{table}



As described in \cite{Rowley3}, lower bounds on $\Gamma(k)$ can derived from the graphs in Table 1 by some simple arithmetic, applying the 2021 construction.  In this way, for example, we show in Table 1 that $\Gamma(5) \ge \sqrt{1428} = 9.919 \dots$ and $\Gamma(6) \ge \sqrt{234} = 15.297 \dots$, and so on.  

An interesting special case also mentioned in \cite{Rowley3} depends on a formula derived in \cite{Abbott} and mentioned in \cite{AbbH} -- namely that $R_r(5) \ge (R_r(3) - 1)^2 + 1$.  A demonstration that $\Gamma(3) \ge 3.280 \dots$ can be found in \cite{Ageron}.  It follows directly that $\Gamma(5) \ge 10.762 \dots$.  

The authors of \cite{XXER} reported that Abbott's result has been extended by Song \cite{Song}.    Song's paper is presented in Chinese, but thankfully, the authors of \cite{XXER} include a proof of the result (along with a further slight extension).

The formula as presented by Song implies that:


\begin{equation}
	R(p_1q_1 + 1, \dots , p_rq_r + 1) > (R(p_1 + 1, \dots , p_r + 1) - 1).(R(q_1 + 1, \dots , q_r + 1) - 1)
	\label{f106}
\end{equation}

This formula allows us to generate lower bounds on $\Gamma(k)$, based on the best available lower bounds for $\Gamma(h)$ for $h < k$.  

However, with the production of the new graphs mentioned in this paper, it is noted that the values of $\Gamma(k)$ shown in Table 1 now generally exceed what is implied by \ref{f106} -- although the bound on $\Gamma(5)$ mentioned above remains a striking special case.


\section{Lower Bounds for $R_3(8)$ and $R_3(9)$}

The existence of specific graphs produced by this research implies new lower bounds for two non-linear `diagonal' 3-colourings of Ramsey graphs.  In each case, these bounds follow from two successive applications of the powerful ``quadrupling'' construction presented as Corollary 5 in \cite{XXER}, in the same manner as in that paper.

Taking the simpler case first, a cyclic non-template (3,8,8; 940)-graph was produced, which can be used to construct a (9,9,9)-graph of order 16 times greater, i.e. equal to 15040, thus proving directly that $R_3(9) \ge 15041$.

The (3,7,7; 673)-graph mentioned in Table 1 can similarly be used to construct a (9,8,8; 10768)-graph.  If the degree of every vertex in the (3,7,7; 673)-graph in the colour 1 is $d_1$, it can be shown that the degree of any vertex in the (9,8,8; 10768)-graph in colour 1 is $9 \times d_1 + 7 \times 673 + 4$.  We observe that the degree of every vertex in the (3,7,7; 673)-graph in colour 1 is 202,  and therefore the degree of every vertex in the (9,8,8; 10768)-graph in colour 1 is 6533.  It follows easily that the vertices in the corresponding neighbourhood of any such vertex induce an (8,8,8)-graph of order 6533.  Hence $R_3(8) \ge 6534$.

The orders of these non-linear graphs are clearly well in excess of what has been possible so far with linear graphs.

\section{Some Conclusions}

The use of SAT-solvers to generate interesting Ramsey graphs adds a new and flexible tool to the search process.  In particular, this approach tends to require less development time than the tree-searching approaches previously used by this author, for any given search.  The outcome in terms of processing times was not so clear.  Graphs with a large number of unconstrained edge-colours produced large input files (to the SAT-solver) -- which in turn required very long processing times, such that achieving a definite result was not certain for large cases.

The specimen graphs found in the past few months result in significant improvements to the lower bounds available for many multicolour (non-diagonal and diagonal) Ramsey numbers, and of course, for the corresponding generalised Schur numbers. 

Cases where $k < 5$ appear least susceptible to these improvements.  The useful template graph with parameters $(3,4,3)$ mentioned above is a notable exception, but the searches have not yet yielded any improvement to the template graphs available with parameters $(4,4,3)$ and $(4,4,4,3)$.  

Much more useful results were achieved in finding $(5,5,5,3)$-template-graphs, and $(k,k,3)$-template-graphs for $5 \le k \le 9$. The latter resulted in improved lower bounds on the orders of larger linear Ramsey graphs $R_r(k)$, exceeding, in most cases, the values shown in the 2021 edition of the Radziszowski Dynamic Survey \cite{RadzDS}.  Given the constraints adopted, it is clear that the searches were designed to be non-exhaustive, but even so the results for $5 \le k \le 9$ represent considerable improvements over previous lower bounds.

In consequence of the identification of these graphs, formula \ref{f106} yields a number of improved lower bounds for $\Gamma(k)$, again for $5 \le k \le 9$, as shown in Table 1.

\subsection*{Acknowledgements}

The production of this paper would have been much more difficult and time-consuming without the ready availability of the SAT-solvers used -- PeneLoPe and Plingeling.  

I express here my sincerest thanks to the developers for making available such powerful tools at no cost. More details of these public-spirited individuals and their creations can be found at \url{http://www.cril.univ-artois.fr/~hoessen/penelope.html} and \url{http://fmv.jku.at/lingeling} respectively.

Once again, I record my warmest thanks to my wife Joan and my son William for their continued support for this work. 



\end{document}